\documentclass[11pt, twoside]{article}

\usepackage{amsmath}    %obligatory package
\usepackage{amssymb}   	%obligatory package
\usepackage{amsthm}     %obligatory package
\usepackage{graphicx}    %obligatory package

\usepackage{textcomp}    %obligatory package
\usepackage[T1]{fontenc} %obligatory package
\usepackage{marvosym}  %obligatory package (for the \Letter symbol for indicating corresponding author)

\usepackage{authblk} 	%obligatory package
\usepackage[usenames]{xcolor} %obligatory package
\usepackage{lastpage} 	%obligatory package

\usepackage[latin1]{inputenc} % allows accented characters									
\usepackage{indentfirst}
\usepackage{amsfonts,amscd}	% symbols and special characters
\usepackage{empheq}
\usepackage{latexsym}		% special characters
\usepackage{graphicx,graphics,epsfig}								
\usepackage{graphpap}										
\usepackage{float}
\usepackage{xcolor}
\usepackage{makeidx}      % remissive index
\usepackage{appendix}	  % appendix
\usepackage{multicol}										
\usepackage{color}

\usepackage[normalem]{ulem}
\usepackage{changepage}
\usepackage{hyperref} % optional arguments: [dvipdfm], [dvips]
\hypersetup{
	colorlinks=true,
	linkcolor=blue,
   citecolor=blue,
   urlcolor=blue
}

\usepackage[dvips=false,pdftex=false,vtex=false]{geometry}
\usepackage{dsfont} % fonts for sets of nymbers, reals, complex, integers, etc.

\numberwithin{equation}{section}
\numberwithin{table}{section}
\numberwithin{figure}{section}

\newcommand{\orcid}[1]{\href{https://orcid.org/#1}{\includegraphics{logo_orcid.jpg}}} % orcid symbol with link to orchid page

\newcommand\raisepunct[1]{\,\mathpunct{\raisebox{0.5ex}{#1}}}
\newtheorem{theorem}{Theorem}[section]

\newtheorem{lemma}[theorem]{Lemma}

\theoremstyle{definition}

\theoremstyle{definition}

\theoremstyle{definition}

\begin{document}

%\maketitle

\thispagestyle{plain}
\renewcommand{\thefootnote}{\arabic{footnote}}   % don't delete!!!
\setcounter{footnote}{0}     % don't delete!!!
\setcounter{page}{1} % The editors will insert the correct initial pagenumber
\label{FirstPage}	

\pagestyle{myheadings} \markboth{\hfil  $\hspace{1.5cm}$  J. L. M Barbosa  $\hspace{1.5cm}$ G. P. Bessa
\hfil $\hspace{3cm}$ } {\hfil$\hspace{1.5cm}$
{Cylinder Theorem}
\hfil}

% Capital letters only for the first letter of the title and for proper names
\title{On the Cylinder Theorem in $M^2\times \mathbb{R}^n $}

\author{\textbf{Jo\~ao L. M. Barbosa \footnote {joaolucasbarbosa@gmail.com}} and {\textbf {G. P. Bessa}}}
%\author[2]{\textbf{Second Author \orcid{222222}}\footnote{The second author is partially supported by.., e-mail: ..}}%Please write the

\affil{Universidade Federal do Ceara, \\ Campus do Pici, Fortaleza, Brazil}

\date{}
\maketitle

%\begin{center}
%\noindent
%\begin{minipage}{0.85\textwidth}\parindent=15.5pt

%\smallskip
%\begin{center}
%\large{\textsl{Dedicated to Professor Renato Tribuzy\\ on the occasion of his 75th birthday}}
%\end{center}

\smallskip

{\small{
\noindent {\bf Abstract.} Consider a surface $M^2$ with Gaussian curvature either $< 0$ or $> 0$. We prove that in $M^2\times \mathbb{R}^n$ cylinders are characterized as the hypersurfaces with both the extrinsic and intrinsic curvatures equal to zero.}
\smallskip

% Please enter at most 6 keywords here with lowercase letters separated by commas.
\smallskip


}
%\end{minipage}

%\end{center}

\section{Introduction} In 1956,  A.V. Pogorelov, \cite{pogorelov} proved that if  $\varphi \colon M \to \mathbb{R}^{3}$ is an isometric immersion of a complete surface with  bounded extrinsic curvature and  locally isometric to the plane then $\varphi(M)$  is a cylinder $\alpha(I) \times \mathbb{R}$, where $\alpha\colon I\subseteq \mathbb{R}$ is a smooth curve. P. Hartman and L. Nirenberg \cite[Thm. III]{HN} extended Pogorelov's result proving that if  $\varphi \colon M \to \mathbb{R}^{n+1}$ is an isometric immersion of a complete $n$-dimensional manifold  with  zero curvature then $\varphi(M)$ is a cylinder $\alpha(I) \times \mathbb{R}^{n-1}$. 
Elementary proofs of the Pogorelov's cylinder theorem   were published independently  by W. S. Massey \cite{Massey} and J. J. Stoker \cite{Stoker}.

In 2020 Barbosa and do Carmo  \cite{lucas2020} proved the cylinder theorem when the ambient space was   $\mathbb{H}^2\times \mathbb{R}$, where $\mathbb{H}^2$ is the $2$-dimension hyperbolic space. This result was extended by J. Park \cite{Park} to include the case when the ambient space was $\mathbb{S}^2\times\mathbb{R}$, where $\mathbb{S}^{2}$ is a sphere.

In the present paper we are extending that result to $M^2\times \mathbb{R}^n$, where $M^2$ is a surface with Gaussian curvature either $< 0$ or $> 0$. In this ambient space we define a {\it cylinder } as the surface given by $\alpha\times \mathbb{R}^n$ where $\alpha$ is a regular curve in $M^2$

We know that any hypersurface in a $N$-dimension space has two curvatures: the intrinsic and the extrinsic Curvature. % that is given by the product of the two principal curvatures, say $k_1k_2$.
In $\mathbb{R}^3$ these two curvatures coincide. In $M^2\times \mathbb{R}^n$ they are in general different and give different information about the hypersurface.

The goal of this paper is to prove the following result
\begin{theorem}\label{MainThm}
 A complete hypersurface  $\Sigma$ in $M^2\times \mathbb{R}^n_{\raisepunct{,}}$ where $M^2$ is a surface with Gaussian curvature either $< 0$ or $> 0$ is a cylinder if and only if  both its intrinsic and extrinsic curvatures vanish.\end{theorem}

%Of course to prove this statement we have to prove two propositions:
%
%\begin{theorem} \label{teo1} Let $\Sigma$ be a complete and connected hypersurface in $M^2\times \mathbb{R}^n$, where $M^2$ is a surface with Gaussian curvature either $< 0$ or $> 0$.  If $\Sigma $ is the cylinder $\alpha\times \mathbb{R}^n$ then it has the intrinsic and extrinsic curvatures equal to zero.
%\end{theorem}
%
%\begin{theorem} \label{teo2} Let $\Sigma$ be a complete and connected hypersurface in $M^2\times \mathbb{R}^n$, where $M^2$ is a surface with Gaussian curvature either $< 0$ or $> 0$.  If $\Sigma$ has both the intrinsic and the extrinsic curvatures equal to zero then $\Sigma$ is a cylinder.
%\end{theorem}

\section{Intrinsic and extrinsic curvatures}

Let $S$ be a Riemannian manifold of dimension $N+2$ with the metric $d\sigma^2$. Let $\Sigma$ be a $N+1$-dimensional submanifold endowed with the induced metric. Let $E_1, \ldots, E_{N+2}$ be an orthonormal frame field in a neighborhood of $S$, adapted to $\Sigma $. Let $\theta_1, \ldots, \theta_{N+2}$ be the corresponding dual forms. The connection forms are defined by
\begin{equation}
d\theta_A = \sum \theta_{AB}\wedge \theta_B, \quad \quad 1\le A,B \le N+2.
\end{equation}
The curvature of $S$ is then defined by
\begin{equation}
{\bar \Omega}_{AB} = d\theta_{AB} - \sum \theta_{AC}\wedge \theta_{CB}
\end{equation}
where $1\le A, B, C \le N+2$.

Restrict to vectors tangent to $\Sigma $ we have
$$
\theta_{N+2} = 0
$$
It follows that
$$
0 = \sum_{B=1}^{N+1}\theta_{N+2\,B} \wedge \theta_B
$$
From this one deduces that
\begin{equation} \label{eq13}
\theta_{N+2\,B} = \sum_{C=1}^{N+1} h_{BC}\theta_C
\end{equation}
where $h_{BC} = h_{CB}$ for $1\le B,C \le N+1$. We say that the direction of $E_A$ is principal when $\theta_{N+2\,A} = h_{AA}\theta_A$. In this case it is common to call $h_{AA}$ principal curvature in the direction of $E_A$ and set $k_A = h_{AA}$.

We also deduce that, for $1\le A,B \le N+1$,
\begin{equation}
{\bar \Omega}_{AB} = d\theta_{AB} - \sum_{C=1}^{N+1}\theta_{AC}\wedge \theta_{CB} - \theta_{A\,N+2}\wedge \theta_{N+2\,B}
\end{equation}
The curvature form in $\Sigma $ is, of course given by
$$
\Omega_{AB} = d\theta_{AB} - \sum_{C=1}^{N+1}\theta_{AC}\wedge \theta_{CB}, \quad \quad \mbox{for } 1\le A,B \le N+1.
$$
The term
\begin{equation} \label{eq1.5}
\Psi_{AB} = \theta_{A\,N+2}\wedge \theta_{N+2\,B}
\end{equation}
is called the { \it extrinsic curvature} of $\Sigma $. We then have
\begin{equation} \label{duascurvaturas}
{\bar \Omega}_{AB} + \Psi_{AB} = \Omega_{AB}, \quad \quad \mbox{ for } 1\le A,B \le N+1.
\end{equation}
If $E_1, \ldots, E_{N+1}$ are principal vectors in $\Sigma$ then $\Psi_{AB} = - k_A\,k_B\theta_A\wedge \theta_B$ and so, equation \ref{duascurvaturas}
becomes
\begin{equation}
{\bar \Omega}_{AB} -k_Ak_B\theta_A\wedge \theta_B = \Omega_{AB}
\end{equation}

\subsection{Intrinsic and extrinsic curvatures of cylinders}\label{secCy}

We are going to compute the intrinsic and the extrinsic curvatures of a cylinder.
Let $\Sigma$ be the cylinder $\alpha \times \mathbb{R}^n$, where $\alpha$ is a regular curve in $M^2$

%To study the geometry of this cylinder we proceed as follows.
%A point in $\alpha \times \mathbb{R}^n$ will be given by $\alpha(s) + (t_1, \ldots, t_n)$.
Observe that the metric in $M^2\times \mathbb{R}^n $ is given by
$$
d\sigma^2 = d\zeta^2 + \sum_{i=1}^n dt_i^2
$$
where $d\zeta $ is the metric in $M^2$ and $\sum dt_i^2$ is the standard metric in $\mathbb{R}^n$.
The covariant differential $\bar D$ in $M^2\times \mathbb{R}^n$ decomposes as
$$
{\bar D} = D + \sum_{i=1}^n d_i
$$
where $D$ is the covariant differential in $M^2$ and the $d_i $'s are  the standard differentials in $\mathbb{R}$.
%if we represent by $\partial /\partial t_i$ the unit vector field tangent to a line in $M^2\times \mathbb{R}^n$, then we obtain
Take an orthonormal frame field given by $e_i = \partial/\partial t_i$, $1\le i \le n$, $e_{n+1} = \alpha'(s)$ and $e_{n+2}$ normal to $\Sigma$, Then
$$
\theta_i = dt_i \quad \quad \mbox{ for } 1\le i \le n, \quad \mbox{ and } \quad \theta_{n+1} = ds
$$
where we assume that $\alpha $ is parameterized by the arc length $s$.
The metric in $\Sigma $ is then given by
$$
\sum_{i=1}^n dt_i^2 + ds^2\,.
$$
Hence $\Sigma $ is isometric to $\mathbb{R}^{n+1} $ and  its intrinsic curvature is zero.

\vspace{0,3 cm}

Given a point $p \in \Sigma$, extend the above mentioned frame field to a neighborhood of $p$ in the ambient space.
Represent by $\omega_{AB}$ the connection forms given by

$$
d\omega_A = \sum \omega_{AB} \wedge \omega_B, \quad \quad 1\le A,B \le n+2.
$$
We then have, for $1\le i \le n$,
$$
0 = d(dt_i) = d\omega_i = \sum_{B=1}^{n+2} \omega_{iB}\wedge\omega_B
$$
Then $\omega_{iB} =0  $,  for $1\le i \le n$ and any $B$.

In particular $\omega_{n+2\, i} = - \omega_{i\, n+2} = 0$ for $1\le i \le n$ and so, from equation (\ref{eq13}),
$$
\sum _{i=1}^{n+1} h_{iB}\omega_B = 0, \quad \quad 1\le B \le n+1\,.
$$
Thus $h_{ij} =0$, in particular $h_{ii} = 0$. Hence, the vectors $e_i$, $1\le i \le n$ are principal vectors and $k_1 = k_2 = \ldots = k_n =0$.
Hence
$$
\Psi_{AB} =0  \quad \mbox{ for } \quad 1\le A,B \le n+1\,.
$$

Therefore,  $\Sigma$ has the intrinsic and the extrinsic curvatures equal to zero.

\section{Proof of theorem \ref{MainThm}}

Let $\Sigma$ be a hypersurface in $M^2\times \mathbb{R}^n$ that we assume to be connected and complete. We also assume that $\Sigma $ is flat and has the extrinsic curvature equal to zero.
Let $p= ({\tilde p}, T)$, ${\tilde p}\in M^2\times \mathbb{R}^n$, ${\tilde p} \in M^2, $ $T = (t_1, \ldots, t_n)$, be a point of $\Sigma $. Take a frame field $e_1, \ldots, e_{n+2}$ adapted to $\Sigma$, that is, restricted to points of $\Sigma $, $e_1, \ldots e_{n+1}$ are vectors tangent to $\Sigma $ and $e_{n+2}$ is normal. Let $\omega_A$, $1\le A \le n+2$, be its dual forms and $\omega_{AB}$, $1\le A,B \le n+2$, be the connection forms associated to such frame.

We then have
$$
d\omega_{AB} = \sum_{C=1}^{n+2} \omega_{AC} \wedge \omega_{CB} + {\bar \Omega}_{AB}
$$
where $\bar \Omega _{A\,B}$ is the curvature form of the plane generated by $e_A, e_B$ in the ambient space.

We now restrict this equation to $\Sigma $ where we have
$$
d\omega_{ij} = \sum_{k=1}^{n+1}\omega_i\wedge \omega_j + \Omega_{ij}, \quad \quad 1\le i,j \le n+1,
$$
where $\Omega$ is the intrinsic curvature of $\Sigma $.

We then have
$$
\Omega_{ij} = {\bar \Omega}_{ij}|_{\Sigma } + \Psi_{ij}, \quad \quad \mbox{ for   } 1\le i,j \le n+1.
$$

Since the intrinsic curvature of $\Sigma $ is zero and the extrinsic curvature of $\Sigma$ is also zero then we conclude that, at the point $p$

$$
{\bar \Omega}_{ij} (e_i, e_j) = 0, \quad \quad \quad 1\le i,j \le n+1
$$

We want to show that $\Sigma$ is of the form $\alpha \times \mathbb{R}^n$. We start by observing that, since the tangent space $T_p(M^2\times \mathbb{R}^n) \approx \mathbb{R}^{n+2}$, then $T_p\Sigma \approx \mathbb{R}^{n+1}$ is either equal to $T_pM^2\times \mathbb{R}^{n-1}$ or the hyperpane $T_p\Sigma $ intersects the tangent plane $T_pM^2$ along a line.

 In the first case the curvature of the plane $T_{\tilde p}M^2 \times \{T\} \subset T_{\tilde p}M^2 \times \mathbb{R}^{n-1}$ is not zero. But, by hypothesis, such curvature should be zero. Therefore this case do not occur.

\vspace{3mm}

To analyze the second case we will prove the following lemma.

\begin{lemma} \label{lemma1}
If $T_p\Sigma $ intersects the plane $T_{\tilde p} M^2 \times \{T\}$ along a line then $\Sigma $ contains ${\tilde p}\times \mathbb{R}^n$.

\end{lemma}

Assume this lemma is proved. As we have seen, for any point $p \in {\Sigma }$ we have that the intrinsic curvature of $T_p\Sigma$ is zero. By what we have seen before, $T_p\Sigma$ will intersect the plane $T_{\tilde p}M^2$ along a line. Now, using the Lemma (\ref{lemma1}) we conclude that $\Sigma $ contains $\mathbb{R}^n$.

Observe that $p\in \Sigma\cap \left(M^2\times \{T\}\right)$. Now, for each $q \in \Sigma\cap M^2\times \{T\}$ there is a one dimensional vecor space $V_q \subset T_qM^2 \times \{T\}$. This gives a one dimensional integrable distribution. Let $\alpha(\mathbb{R})\subset M^2\times \{T\}$ be the maximal integrable curve of this distribution. By the Lemma (\ref{lemma1}), $\alpha(t) \times \mathbb{R}^n \subset \Sigma$ for every $t\in \mathbb{R}$. This shows that $\alpha(\mathbb{R}) \times \mathbb{R}^n \subset \Sigma $. Since $\Sigma $ is a hypersurface we have $\Sigma = \alpha(\mathbb{R}) \times \mathbb{R}^n$.

%Choose constant vector $V\in \mathbb{R}^n$. Since the choice of $V$ is possible at each point $p \in \Sigma$, the vector $V$ %defines
%a vector field in $\Sigma$ whose integral curves are lines. Such lines are, of course, contained in $\Sigma$. Since this is true %for any choice of $V$ one concludes that $\mathbb{R}^n \subset \Sigma$. It then follows that $\Sigma = \alpha \times  %\mathbb{R}^n$. Since $\Sigma $ is a hypersurface then one concludes that $\alpha $ is regular.

Therefore $\Sigma $ is a cylinder. This proves Theorem (\ref{MainThm})) since cylinder has intrinsic and extrinsic curvatures as shown in section \ref{secCy}.

\vspace{3 mm}

\subsection{Proof of the Lemma \ref{lemma1}}
Choose in neighborhood of $p=(\tilde{p},T) \in M^2\times \mathbb{R}^n_{\raisepunct{,}}$  a frame field  given by $E_1 = \partial/\partial t_1, \ldots, E_n = \partial/\partial t_n$, and $E_{n+1}$, $E_{n+2}$ tangent to $M^2\times \{T\}$, $T= (t_1, \ldots, t_n)$. Let $\omega_i = dt_i$, $1\le i \le n$, and $\omega_{n+1}$, and $\omega_{n+2}$ the dual forms of $E_{n+1}$, $E_{n+2}$. Let $\omega_{AB}$ be the corresponding connection forms defined by
$$
{\bar D}E_A = \sum \omega_{AB}E_B, \quad \quad 1\le A,B \le n+2.
$$
Observe that,

$$
\omega_{iA}=0, \quad \quad 1\le i \le n \quad 1\le A \le n+2.
$$
 It then follows that
\begin{equation} \label{eq1323}
0 = d\omega_{iA} = {\bar \Omega} _{iA}, \quad \quad 1\leq i \leq n, \quad 1\le A \le n+2.
\end{equation}
%It also implies that
%\begin{equation} \label{eq1323}
%{\bar \Omega}_{i\,n+1} = {\bar \Omega}_{i\,n+2} =0, \quad \quad 1\le i \le n.
%\end{equation}

Set $c(p)= {\bar \Omega}_{n+1\,n+2}(E_{n+1},E_{n+2})$ at the point $p$. Then
\begin{equation}
{\bar \Omega}_{n+1\, n+2} = c(p)\,\omega_{n+1}\wedge \omega_{n+2}
\end{equation}

\vspace{3mm}

We are interested to compute the curvature of any $2$-plane  of $T_p\Sigma $ at some point $p = ({\tilde p},T)$, $T = (t_1, \ldots, t_n)$. If $T_p\Sigma $ is vertical, that is, if $  \{\tilde p\} \times\mathbb{R}^n\subset T_p\Sigma$, we choose an orthonormal basis $V_i = \partial/\partial t_i $, $1\le i \le n$, $V_{n+1}$ of $
T_p\Sigma $. Then $V_i = E_i$, $1\le i \le n$ and $V_{n+1} = a\,E_{n+1} + b\,E_{n+2}$, $a^2 + b^2 = 1$. From (\ref{eq1323}) it follows that
$$
\Omega_{A\,B}(V_A, V_B) = 0 \quad \quad 1\le A,B \le n+1.
$$

If $T_p\Sigma $ is not vertical we may have two distinct situations:

{\bf Situation A}
$$
T_p\Sigma = {\cal P}\times \mathbb{R}^{n-1}
$$
where ${\cal P}$ is a plane tangent to $M^2\times \{T\}$, $T=(t_1,\ldots,t_n)$, that obviously has the same curvature of $T_{p}(M^2\times \{T\})$, whose value is different from zero. But, by hypothesis, the curvature of any plane in $T\Sigma$  should be zero. Therefore this case do not occurs.

\vspace{3mm}

{\bf Situation B}

We also have
$$
T_p\Sigma = {\cal P}\times \mathbb{R}^{n-1}
$$
where ${\cal P}$ is a plane that cuts $T_p(M^2\times \{T\})$ along a line generated by a unit vector $V\in T_p(M^2\times \{T\})$. Observe that $V$ is perpendicular to $\{\tilde p\}\times \mathbb{R}^n$. In this case we may take a frame field $E_1,\ldots,E_{n+2}$ in a neighborhood of $p$ choosing
$E_i = \partial/\partial t_i$, $1\le i \le n$, $E_{n+2} = V$, $E_{n+1}$ normal to all vectors chosen.

Take another frame field $e_1,\ldots,e_{n+2}$ this time given by $e_i = \partial/\partial t_i $, $1\le i \le n-1$, $e_{n+2} = V$, $e_n$ and $e_{n+1}$ in the plane generated by $E_n$ and $E_{n+1}$. This choice made such that, in $p$, $e_n$ is in $T_p\Sigma $ and $e_{n+1}$ is normal to $T_p\Sigma $. Thus $T_p\Sigma $ is generated by $e_1, \ldots, e_n, e_{n+2}$ and has $e_{n+1}$ as a normal vector.
In this way we must have
\begin{eqnarray} \nonumber
e_n & = & a\,E_n + b\,E_{n+1} \\
e_{n+1} & = & -b\,E_n + a\,E_{n+1}
\end{eqnarray}
being $a^2 + b^2 = 1$.
Of course we also have
\begin{eqnarray} \nonumber
E_n &=& a\,e_n - b\,e_{n+1} \\
E_{n+1} &=& b\,e_n + a\,e_{n+1}
\end{eqnarray}
Let $\theta_A$, $1\le A \le n+2$, be the dual forms of the frame field $e_A$, while $\omega_A$ are the dual forms of the frame field $E_A$. By the choices made we have:
\begin{eqnarray} \nonumber
\omega_i &=& dt_i \quad \quad 1\le i \le n \\ \nonumber
\theta_i &=& dt_i \quad \quad 1\le i \le n-1\,.
\end{eqnarray}
Thus, from the first equation,  we have ${\bar D}E_i =0$ and so $\omega_{iA} =0$, $1\le i \le n$ and any $A$.  From the second one
we have $\theta_{iA} = 0$ for $1\le i \le n-1$ and any $A$. Therefore, it can be different from zero only the connection forms
$\omega_{n+1\,n+2}$, $\theta_{n\,n+1}$, $\theta_{n,n+2}$ and $\theta_{n+1\,n+2}$.

Since $e_{n+2} = E_{n+2}$ then we have
\begin{eqnarray} \nonumber
{\bar D}e_{n+2} &=& {\bar D}E_{n+2} = \omega_{n+2\,\, n+1}E_{n+1} \\
&=& \omega_{n+2\,\,n+1}(b\,e_n + a\,e_{n+1}) \\ \nonumber
&= &\theta_{n+2\,\,n}e_n + \theta_{n+2\,\,n+1}e_{n+1}
\end{eqnarray}
From this it follows that
\begin{equation} \label{res01}
\theta_{n+2\,\,n} = b\,\,\omega_{n+2\,\,n+1} \quad \quad \theta_{n+2\,\,n+1} = a\,\,\omega_{n+2\,\,n+1}
\end{equation}
We compute
\begin{eqnarray} \nonumber
{\bar D}e_{n+1} &=& {\bar D}(-bE_n + aE_{n+1})\\ \nonumber & =& -b{\bar D}E_n + a{\bar D}E_{n+1} \\ \nonumber
&=& a(\omega_{n+1\,\,n}\,E_n + \omega_{n+1\,\,n+2}\,E_{n+2}) \\
&=& a\omega_{n+1\,n+2}E_{n+2} = a\omega_{n+1\,n+2}e_{n+2} \\ \nonumber
&=& \theta_{n+1\,n}\, e_n + \theta_{n+1\, n+2} \,e_{n+2}
\end{eqnarray}
From this we deduce that
\begin{equation} \label{res02}
\theta_{n+1\,\,n} = 0 \quad \quad \theta_{n+1\,\,n+2} = a\,\omega_{n+1\,\,n+2}
\end{equation}
Similarly we have
\begin{eqnarray} \nonumber
{\bar D}e_n &=& {\bar D}(a\,E_n + b\,E_{n+1}) \\ \nonumber & = & a\,{\bar D}E_n + b\,{\bar D}E_{n+1}  \\ \nonumber
&=& b \,(\omega_{n+1\,\,n}E_n + \omega_{n+1\,\,n+2}E_{n+2})  \\
&=& b\,\omega_{n+1\,n+2}\,E_{n+2} = b\,\omega_{n+1\,n+2}\,e_{n+2}  \\ \nonumber
&=& \theta_{n\,\,n+1} \,e_{n+1} + \theta_{n\,\, n+2} \,e_{n+2}
\end{eqnarray}
From this it follows that
\begin{equation} \label{res03}
\theta_{n\,\,n+1} = 0 \quad \quad \theta_{n\,\,n+2} = b\,\omega_{n+1\,\,n+2}
\end{equation}

We will represent by $R$ the curvature tensor of $M^2\times \mathbb{R}^n_{\raisepunct{,}}$ so that
$$
R(e_i,e_j, e_m, e_k) = {\bar \Omega}_{ij}(e_m,e_k).
$$
We now compute the value of the curvature at $T_p\Sigma $, that, remember, is generated by $e_1, \ldots, e_n, e_{n+2}$. Using the $4$-linearity of $R$, we obtain, for $1\le i,j \le n-1$,
$$
R(e_i, e_A, e_i, e_A) = R(E_i,E_A,E_i,E_A) = {\bar \Omega}_{i\,A}(E_i,E_A) = 0
$$
were we have used $\omega_{i\,A} = 0$ for $1\le i \le n-1$, $1\le A \le n+2$. Hence, the only possible non-zero cases, in the computation of the curvature, are
$$
R(e_n, e_{n+1}, e_n, e_{n+1}), \quad R(e_n, e_{n+2}, e_n, e_{n+2}), \quad R(e_{n+1}, e_{n+2}, e_{n+1}, e_{n+2}).
$$
\begin{eqnarray} \nonumber  \label{eq312}
 R(e_n,e_{n+2},e_n,e_{n+2})& = & R(a\,E_n + b\,E_{n+1}, E_{n+2}, a\,E_n + b\,E_{n+1}, E_{n+2}) \\  \nonumber
& = & a^2\,R(E_n,E_{n+2},E_n,E_{n+2})\nonumber \\ && +\, ab\,R(E_n, E_{n+2},E_{n+1},E_{n+2})  \\ \nonumber
&&+ \, ba\,R(E_{n+1},E_{n+2},E_n,E_{n+2})\nonumber \\ && +\, b^2 \,R(E_{n+1},E_{n+2},E_{n+1},E_{n+2})  \nonumber
\end{eqnarray}
Let $c$ be the value of the curvature of $M^2\times \{T\}$ at the point $p= ({\tilde p}, T)$.
Observe that
\begin{eqnarray} \nonumber
R(E_n,E_{n+2},E_n,E_{n+2}) &= &{\bar \Omega}_{n\,\,\,n+2}\,(E_n,E_{n+2}) = c, \\ \nonumber
R(E_n,E_{n+2},E_{n+1},E_{n+2}) &= &{\bar \Omega}_{n\,\, n+2}\,(E_{n+1},E_{n+2}) = 0, \\ \nonumber
R(E_{n+1},E_{n+2},E_n,E_{n+2}) &= & R(E_n,E_{n+2},E_{n+1},E_{n+2}) = 0 \\ \nonumber
\end{eqnarray}
and
\begin{eqnarray} \nonumber
R(E_{n+1},E_{n+2}, E_{n+1}, E_{n+2})& = & {\bar \Omega}_{n+1\, n+2}(E_{n+1},E_{n+2}) = 0 \\ \nonumber
\end{eqnarray}
In the last equality we have used (\ref{eq1323}). Thus, from equation (\ref{eq312}) we obtain

\begin{equation}
R(e_n,e_{n+2},e_n,e_{n+2}) = a^2c
\end{equation}

Since $c\neq 0$ then $R(e_n,e_{n+2},e_n,e_{n+2}) = 0$ if and only if $a=0$. But this is equivalent to
the fact that the vector $\partial/\partial t_n$ belongs to $T_p\Sigma $, besides $\partial/\partial t_1, \ldots , \partial/\partial t_{n-1}$. So, $T_p\Sigma \supset \{\tilde p \} \times \mathbb{R}^n $.  Hence, $\Sigma $ contains $\mathbb{R}^n$.
This concludes de proof of the Lemma.

%\strictpagecheck
%\checkoddpage
%\ifoddpage
%\newpage{\ }\thispagestyle{empty}
%\fi

\end{document}